\documentclass[smallcondensed]{svjour3}     
\smartqed  
\usepackage[caption=false]{subfig}
\usepackage[utf8]{inputenc}
\usepackage{graphicx}
\usepackage{txfonts}%
\usepackage{color}
\begin{document}

\title{A Lower Dimensional Linear Equation Approach to The M-Tensor Complementarity
Problem\thanks{Supported by the NSF of China grant 11771157 and the NSF of Guangdong Province
grant No.2020B1515310013.}}


\author{Dong-Hui Li \and Cui-Dan Chen \and Hong-Bo Guan}

\institute{Dong-Hui Li, \at
              School of Mathematical Sciences, South China Normal University, Guangzhou, 510631, Guangdong, China.\\
              \email{lidonghui@m.scnu.edu.cn}
              \and
         Cui-Dan Chen \at
              School of Mathematical Sciences, South China Normal University, Guangzhou, 510631, Guangdong, China.\\
              \email{Cuidan.Chen@m.scnu.edu.cn}
              \and
          Hong-Bo Guan, Corresponding author\at
              $^1$ School of Mathematical Sciences, South China Normal University, Guangzhou, 510631, Guangdong, China.\\
              $^2$ School of Mathematics, Physics and Energy Engineering, Hunan Institute of Technology, HengYang, 421002, Hunan, China. \\
              \email{hongbo\_guan@m.scnu.edu.cn}
}
\date{Received: date / Accepted: date}

\maketitle

\begin{abstract}
We are interested in finding a solution to the tensor complementarity problem  with a strong M-tensor,
which we call the M-tensor complementarity problem.
We propose a lower dimensional linear equation approach to solve that problem.
At each iteration, only a lower dimensional system of linear equation needs to be solved.
The coefficient matrices of the lower dimensional linear systems
are independent of the iteration after finitely many iterations.
We show that starting from zero or some nonnegative point,
the method generates a sequence of iterates that converges to a  solution of the
problem monotonically. We then make an improvement to the method and establish its monotone convergence.
At last,
we do numerical experiments to test the proposed methods. The results positively support
 the proposed methods.
\keywords{ M-tensor complementarity problem \and lower dimensional linear equation approach \and  monotone convergence}
\subclass{15A69 \and 65K10 \and 90C33}
\end{abstract}

\section{Introduction}
\label{intro}
As a direct extension of the linear complementarity problem (LCP),
the tensor complementarity problem (TCP) is a subclass of the general nonlinear complementarity
problem (NCP) with some attractive special structures. In recent years, the study in the TCP has
received much attention due to its wide applications in multi-person noncooperative games,
hypergraph clustering problem, DNA microarrays  and so on, see, for instance,
\cite{Bulo-Pelillo-09,Huang-Qi-17,Luo-Qi-Xiu-17} and references therein.
We refer to \cite{Huang-Qi-JOTA-1,Huang-Qi-JOTA-3,Qi-Chen-Chen-18} for good reviews in
the theory and the applications of  TCPs.

In this paper, we focus our attention on numerical methods for solving TCPs.
Since  TCPs are  NCPs, in many cases, they can be
solved by existing efficient numerical methods for solving NCPs.
Indeed, those methods such as the smoothing Newton method \cite{Huang-Qi-17}, the nonsmooth Newton method \cite {Liu-Li-Vong-18},
the Levenburge-Marquardt method \cite{Du-Zhang-Chen-Qi-18},
continuation method \cite{Han-19}  and the interior point method \cite{Zhang-Chen-Zhao-19} for solving NCP have been applied
to solve TCPs.
Another way to develop numerical methods for solving TCPs is to
extend the idea of those methods for solving LCPs. For example,
Xie, Li and Xu \cite{Xie-Li-Xu-17} proposed an iterative method for finding the least solution to the TCP
with a strongly monotone Z-tensor where the subproblems are lower dimensional
tensor equations. After finitely many iterations, the solution of the last
tensor equation is a solution of the TCP. Similar to LCP, the TCP can be formulated as
a  modulus-based nonsmooth equation and then  the idea of numerical methods for
solving the nonsmooth equation reformulation to LCP can be applied \cite{Liu-Li-Vong-18,Xu-Li-Xie-19}.
Recently, Du and Zhang \cite{Du-Zhang-19} showed that a TCP can be transformed as a mixed integer programming.
There are also some other interesting ideas for developing numerical methods to TCPs
\cite{Xie-Jin-Wei-2018a,Zhao-Fan-18}.
We refer to a recent review paper \cite{Huang-Qi-JOTA-2} for a good summary in the numerical methods
for  TCPs.

Quite recently, Guan and Li \cite{Guan-Li-19} proposed a linearized method and a lower dimensional
linear equation method for solving the TCP with a strong M-tensor (M-TCP).
Monotone convergence for both methods was well established.
However, to ensure the convergence of the methods, the initial point must be selected in
some feasible set.
It really restricts the application to the methods because in many case, finding
a feasible point is not an easy task. In this paper, we
 further study numerical methods for finding a solution to the M-TCP.
The proposed method can be regarded as an improvement to the method in \cite{Guan-Li-19}. The  subproblems of the method are
also lower
dimensional systems of linear equations.
The proposed method possesses the following attractive properties.

\begin{itemize}
\item At each iteration, only a lower dimensional system of linear equations needs to be solved.
\item The initial point can be set to zero or  a solution of some lower dimensional
system of linear equations.
\item The sequence of the generated iterates  converges to a solution of the M-TCP monotonically.
\end{itemize}

We then make an improvement to the method so that the improved
method is closer to Newton's method than its unimproved version.
 We prove that the improved method still possesses monotone convergence property.

The structure of this paper is organized as follows. In the next section,
we first introduce some necessary concepts and notations and then propose
a lower dimensional linear equation method.
We will also establish the monotone convergence of the method. In Section \ref{sec:3},
we make an improvement
to the method proposed in Section 2 and establish its monotone convergence.
In Section \ref{sec:4}, we conduct some numerical experiments to test the effectiveness
of the proposed methods.

We conclude this section by introducing some notations that will be used throughout the paper.
Let $\mathbb{R} ^n$ denote the $n$-dimensional real Euclidean space.
Denote $[n]=\{1,2, ... ,n\}$.
For a column vector $x=(x_1,x_2,...,x_n)^T\in \mathbb R^{n}$ and $I\subseteq[n]$,
$x_I$ denotes the subvector of $x$ whose elements are  $x_i\in \mathbb R$, $i\in I$.
For matrix $A=(a_{ij})\in  \mathbb R^{n\times n}$ and $I,J\subseteq [n]$,
we denote $A_{IJ}$ the submatrix of $A$ with elements $a_{ij}$, $i\in I,j\in J$.
For convenience, the principle submatrix $A_{II}$ is abbreviated as $A_{I}$.
Similarly, for a tensor ${\cal A}=(a_{i_1i_2...i_m})$ and index set $I\subset [n]$,
we use ${\cal A}_I$ to denote the subtensor of ${\cal A}$
with elements $a_{i_1i_2...i_m}$, $i_1,i_2,...,i_m\in I$.

\section{A Lower Dimensional Linear Equation Approach To M-TCP}\label{sec:1}

We first introduce the some concepts which will be used in this paper \cite{Ding-Qi-Wei-13,Qi-Chen-Chen-18,Qi-Luo-17}.

\begin{definition}\label{tensor}
A real tensor (hypermatrix) ${\cal A} = (a_{i_1\ldots i_m} )$ is a multi-array of real entries $a_{i_1\ldots i_m}\in \mathbb {R}$,
where $i_j = 1, \ldots, n_j$ for $j = 1, \ldots,m$. When $n_1 = \cdots = n_m = n$,
$\cal A$  is called an $m$th order $n$-dimensional tensor.
The set of all $m$th order $n$-dimensional
real tensors is denoted as $ {\cal T}(m,n)$.
\end{definition}

Clearly, when $m=2$, ${\cal T}(2,n)$ is the set of all $n\times n$ real matrices.
For a tensor ${\cal A}\in {\cal T}(m,n)$ and a vector $x\in \mathbb {R}^n$, the notation ${\cal A}x^{m-1}$ stands
for a vector in $\mathbb {R}^n$ whose elements are
\[
\sum\limits_{i_2,\ldots, i_m}a_{ii_2\ldots i_m}x_{i_2}\cdots x_{i_m},\quad i\in [n],
\]
 and
\[
{\cal A}x^m=x^T\Big ({\cal A}x^{m-1} \Big )=\sum\limits_{i_1i_2,\ldots, i_m}a_{i_1\ldots i_m}x_{i_1}x_{i_2}\cdots x_{i_m}.
\]

\begin{definition}\label{tens:P}
Let  ${\cal A}=(a_{i_1i_2,...,i_m})\in {\cal T}(m,n)$.
\begin{itemize}
\item Tensor $\cal A$ is called nonnegative and denoted as ${\cal A}\ge 0$ if all its elements are nonnegative.
\item If the elements of $\cal A$ satisfy $a_{i\ldots i}=1$, $i\in [n]$, and all other elements are zeros,
then $\cal A$ is called the identity tensor and is denoted by $\cal I$.
\item Tensor $\cal A$ is called a P-tensor if for any vector $x\in \mathbb {R}^n\setminus\{0\}$,
\[
\max_{i\in[n]}x_i({\cal A}x^{m-1})_i>0.
\]
\item Tensor $\cal A$ is called a Z-tensor if all of its
off-diagonal entries are nonpositive, which is equivalent to
\begin{equation}\label{tens:M-Z}
{\cal A} = s{\cal I} - {\cal B},\quad {\cal B}\geq 0,
\end{equation}
where $s>0$ is a scalar.
In the case $m=2$, Z-tensors reduce to the so-called Z-matrices \cite{Abraham-94}.
\item A Z-tensor of form (\ref{tens:M-Z}) is called an M-tensor if $s\geq \rho({\cal B})$.
 If $s> \rho({\cal B})$, then ${\cal A}$ is called a strong or nonsingular M-tensor.
 Here and throughout the paper, $\rho({\cal B})$ denotes the the spectral of tensor $\cal B$ \cite{Qi-Luo-17,Qi-Chen-Chen-18}.
 In the case $m=2$, M-tensors reduce to the so-called M-matrices \cite{Abraham-94}.
\end{itemize}
\end{definition}

The TCP was first introduced by Song and Qi \cite{Song-Qi-16}. It is to find an $x\in\mathbb  R^n$ such that
\begin{equation}\label{M-tcp}
x\ge 0,\quad F(x)={\cal A}x^{m-1}-b\ge 0,\quad x^TF(x)=0.
\end{equation}
where ${\cal A}\in{\cal T}(m,n)$ and $b\in \mathbb {R}^n$.
If ${\cal A}$ is a strong M-tensor, then the problem is an M-TCP. Since  strong M tensors are
 P tensors, the M-TCP (\ref{M-tcp}) always  has a solution \cite{Qi-Chen-Chen-18}.
 In particular, the so called feasible set
 \[
 {\cal F}=\{x\in \mathbb R^n\;|\; F(x)\ge 0,\; x\ge 0\}
 \]
has a least element that is a solution to the M-TCP.

\begin{definition}
For a tensor ${\cal A}=(a_{i_1\ldots i_m})\in {\cal T}(m,n)$,
we call the subtensor ${\cal M}({\cal A})=(\tilde a_{i_1\ldots i_m})\in {\cal T}(m,n)$
the majorization tensor of $\cal A$
whose elements are
\[
\tilde a_{ij\ldots j}=a_{ij\ldots j},\qquad i,j=1,2,\ldots,n
\]
and all other elements are zeros.
The corresponding matrix $M({\cal A})=(a_{ij})\in \mathbb R^{n\times n}$ with
\[
a_{ij}=a_{ij\ldots j},\qquad i,j=1,2,\ldots,n
\]
is called the majorization matrix of $\cal A$. The subtensor
\[
{\cal M}_c({\cal A})={\cal A}-{\cal M}({\cal A})
\]
is called the complement of $\cal A$ with respective to ${\cal M}({\cal A})$.
\end{definition}

For any $x=(x_1,\ldots,x_n)^T\in\mathbb  R^n$ and $\alpha \in \mathbb R$, we
define $x^{[\alpha]}=(x_1^{\alpha},\ldots, x_n^\alpha)^T$
whenever it is meaningful.
It is easy to get
\[
{\cal A}x^{m-1}= M({\cal A}) x^{[m-1]}+ {\cal M}_c({\cal A}) x^{m-1}.
\]
For simplicity, in the latter parts of the paper, without confliction, we simply denote the majorization matrix
$ M({\cal A})$ by $A$ and the subtensor ${\cal M}_c({\cal A})$ by $\bar {\cal A}$, i.e.,
\[
A= (a_{ij})=M({\cal A})=(a_{ij\ldots j}),\quad  \bar {\cal A}= {\cal M}_c({\cal A}).
\]
By the use of those notations, the TCP (\ref{M-tcp})  can be written as
\begin{equation}\label{eqn:M-tcp-EQU}
x\ge 0,\quad F(x)=A x^{[m-1]}+ \bar {\cal A} x^{m-1}-b\ge 0,\quad x^TF(x)=0.
\end{equation}

We are going to propose a sequential lower dimensional linear equation method to solve the last problem.
For an $x\in\mathbb  R^{n}_{+}=\{x\in\mathbb  R^n|\;x\geq 0\}$ , we denote
\[
\tilde I_+(x)=\{i\ \ |x_i>0\} \quad \mbox{and}\quad  \tilde I_0(x)=\{i\ \ |x_i=0\}.
\]
It is not difficult to see that if we have known a solution $x^{*}$ of the
TCP (\ref{M-tcp}), then  $x^{*}$ is a solution of the following lower dimensional tensor equation
\[
F_i(x)=0,\ \ \forall i\in \tilde I_+(x^{*}), \quad x_i = 0,\ \ \forall i\in \tilde I_0(x^{*}).
\]
One the other hand, if for some index set $I$, the lower dimensional tensor equation
\begin{equation}\label{equiv:NCP}
F_i(x)=0,\ \ \forall i\in I, \quad x_i=0,\ \ \forall i\notin I
\end{equation}
has a nonnegative solution $\bar{x}$ satisfying
\[
F_i(\bar{x})\geq 0,\ \ \forall i\notin I,
\]
then $\bar{x}$ is a solution to the TCP.

The above equivalency between the NCP/LCP and the system lower dimensional nonlinear/linear equations
provides a way to develop lower dimensional equation approach to the NCP/LCP. Indeed,
the LCP with a Z-matrix can be solved via solving several systems of linear equations (see
e.g. \cite {Li-Nie-Zeng-08}). Such an idea has recently extended to solving the Z-tensor complementarity
problem (Z-TCP) by Xie, Li and Xu \cite{Xie-Li-Xu-17} where the Z-TCP was solved by solving
several lower dimensional tensor equations. Our purpose here is to improve the method by
Xie, Li and Xu \cite{Xie-Li-Xu-17}. Specifically, we will propose an iterative method
to  solve the M-TCP by solving a sequence of systems of linear equations.

From now on, without specification, we always assume that $\cal A$ is a strong M-tensor.
It is easy to see that if $\cal A$ is a strong M-tensor, then its majorization matrix $A={\cal M}({\cal A})$ is a nonsingular
M-matrix, in particular, $A^{-1}$ exists and is nonnegative.

The following lemma proved by Li, Guan and Wang \cite{Li-Guan-Wang-18} will be very important role
in the development of our method.
\begin{lemma}\label{lem-g}
If ${\cal A}$ is a strong M-tensor, and the feasible set ${\cal S}$ defined by
\[
{\cal S} \stackrel\triangle {=}\{x\in \mathbb R^n_+|\;F(x)={\cal A}x^{m-1}-b\leq 0\}
\]
is not empty, then ${\cal S}$ has a largest element that is the largest
nonnegative solution to the M-tensor equation $F(x)={\cal A}x^{m-1}-b=0$.
\end{lemma}

We are in the position to describe the idea  of our method, which we call
sequential lower dimensional linear equation method.
At the beginning, we get an initial point $x^{(0)}\ge 0$ as an estimate to $x^*$ and an
initial index set $I_0=\{i\;|\; F_i(x^{(0)})<0\}\neq \emptyset$  as an estimate to $\tilde I_+(x^{*})$ satisfying
$I_0\subseteq \tilde I_+(x^*)$.
In general, at iteration $k$, suppose we have got an $x^{(k)}$ such that
$I_k=\{i\;|\; F_i(x^{(k)})<0\}\neq \emptyset$.
We then solve a lower dimensional system of linear equations
to get the next iterates $x_{I_k}^{(k+1)}$ and let
$x_i^{(k+1)}=x_i^{(k)}=0$, $\forall i\not\in I_k$.
As we shall show that the sequences of iterates $\{x^{(k)}\}$ and index sets $\{I_k\}$
will satisfy
\[
0\le x^{(k)}\le x^{(k+1)}\quad\mbox{and}\quad I_k\subseteq I_{k+1}.
\]
After finitely many iterations, the index set $I_k$ will remain unchanged.
That is, there is an nonnegative integer $\bar k$ such that
$I_k=I_{\bar k}$, $\forall k\ge \bar k$. As a result,
the method essentially reduces to a linearized method for finding a
nonnegative solution to the lower dimensional nonlinear equation
\[
F_i(x)=0,\quad i\in I_{\bar k}.
\]
Since $x_i^{(k)}=0$ and $F_i(x^{(k)})\ge 0$, $\forall i\not\in I_{\bar k}$,
any nonnegative solution of the last lower dimensional equation extended with
some zero elements  is a solution of the M-TCP (\ref{M-tcp}).

The steps of the method are given below.

\noindent
{\bf{Algorithm 1. (A Sequential Lower Dimensional Linear Equation Approach)}}
\begin{itemize}
\item [] {\bf Initial.} Given positive sequence $\{\alpha_k\}\subset (\alpha_{\min},1)$
and positive constant $\eta > 0 $. Find an initial point  $x^{(0)}\ge 0$ such that
the set $I_0=\{i\;|\;F_i(x^{(0)})<0\}$ is not empty. Let $k:=0$.
\item [] {\bf Step 1.} Stop if  $\|\min\{F(x^{(k)}),x^{(k)}\}\|\leq \eta$.
 Otherwise,
 solve the following lower dimensional system of linear equations
\begin{equation}\label{sub:Meq}
A_{I_k} (x_ {I_k}^{(k+1)})^{[m-1]}
=A_{I_k}\Big (x_{I_k}^{(k)}\Big ) ^{[m-1]}-\alpha_k F_{I_k}(x^{(k)})
\end{equation}
 to get $x^{(k+1)}_{I_k}$.
 Let $x_i^{(k+1)}=x_i^{(k)}=0$, $\forall i\not\in I_k$.
\item [] {\bf Step 2.}  Determine  the index set
\[
I_{k+1}=\{i\;|\; F_i(x^{(k+1)})< 0\}.
\]
Let $k:=k+1$. Go to Step 1.
\end{itemize}

\noindent
{\bf Remark 1.}
At the beginning, we need to find an $x^{(0)}\ge 0$ such that $I_0\neq\emptyset$
and $I_0\subseteq \tilde I_+(x^*)$, where $x^*$ is a solution to the M-TCP.
To this end, we define
the index set
\[
 I_0=\{i\;|\; b_i>0\}.
\]
If $I_0=\emptyset$, then zero is a trivial solution to the M-TCP (\ref{M-tcp}).
Without loss of generality, we suppose that there are at least one
index $i$ satisfying $b_i>0$, namely, $I_0\neq\emptyset$.
It is easy to see that  the relation $ I_0\subseteq \tilde I_+(x^*)$ holds for any solution $x^*$
of (\ref{M-tcp}).
In this way, we can easily get an initial point $x^{(0)}=0$ and $I_0$.

We may also consider to find a larger initial point by solving a lower dimensional
  system of linear equations
\[
A_{I_0}x_{I_0}^{[m-1]}-b_{I_0}=0
\]
to get $x^{(0)}_{I_0}$ and then let $x^{(0)}_i=0$, $\forall i\not\in I_0$.
It is not difficult to see that such an $x^{(0)}$ meets the requirement
that $x^{(0)}\ge 0$ and $\{i\;|\;F_i(x^{(0)})<0\}\neq \emptyset$ unless
\[
a_{i_1i_2\ldots i_m}=0,\quad \forall (i_1,i_2,\ldots,i_m)\neq (i,j,\ldots,j),\, i,j\in I_0.
\]
In the latter case, we can select $\rho x^{(0)}$ with $\rho \in (0,1)$ as an
initial point.

In what follows, we investigate some interesting properties of the above algorithm.

Since $A$ is a nonsingular M-matrix, so is its principal submatrix $A_{I_k}$.
Therefore, the system of linear equations (\ref{sub:Meq}) always has a unique solution
unless $I_k=\emptyset$. As we shall show in Proposition \ref{prop:monotone-1} that
the index set $I_k$ will never be empty. So, the algorithm is well defined.

By the fact $F_{I_k}(x^{(k)})<0$, we can easily get the following trivial proposition.
\begin{proposition}\label{prop:monotone-1-c}
Suppose that ${\cal A}$ is a strong M-tensor. Then the sequence $\{x^{(k)}\}$ generated
by the sequential lower dimensional linear equation approach  is non-decreasing in the sense
\begin{equation}\label{monotone-1-c}
x^{(k+1)}\ge x^{(k)}\ge 0,\quad\forall k\ge 0.
\end{equation}
\end{proposition}

The proposition below show  that
the sequence of index sets $\{I_k\}$ is non-decreasing, i.e.,
$I_k\subseteq I_{k+1}$.
As a result, the method is well defined. Since $I_k\subseteq [n]$ is finite, after finitely many iterations, the
index set $I_k$ will remain unchanged.

\begin{proposition}\label{prop:monotone-1}
Suppose that ${\cal A}$ is a strong M-tensor. Then the sequence of the  index sets $\{I_k\}$ generated by  Algorithm 1
is non-decreasing, i.e.,
\begin{equation}\label{Ik-c}
I_k\subseteq I_{k+1}.
\end{equation}
In particular, there is an index $\bar k\ge 0$ such that
\[
I_k = I_{\bar k}\stackrel\triangle {=} \bar I.
\]
Moreover, it holds that
\begin{equation}\label{dec-F}
 F_i(x^{(k+1)})\le F_i(x^{(k)}),\quad\forall i\not\in I_k.
\end{equation}
\end{proposition}
{\bf{Proof}}
We first verify (\ref{Ik-c}). Indeed, we have
\begin{eqnarray*}
F_{I_k}(x^{(k+1)}) &=& \Big ({\cal A} \Big (x^{(k+1)}\Big )^{m-1}-b\Big ) _{I_k} \\
    &=&  A_{I_k}\Big (x_{I_k}^{(k+1)}\Big )^{[m-1]}+\Big ( \bar {\cal A} \Big (x^{(k+1)}\Big )^{m-1}\Big )_{I_k}-b_{I_k}\\
    &=&   A_{I_k}\Big (x_{I_k}^{(k)}\Big )^{[m-1]}-\alpha _kF_{I_k}(x^{(k)})
                +\Big ( \bar {\cal A}\Big (x^{(k+1)}\Big )^{m-1}\Big )_{I_k} -b_{I_k}\\
    &=&  (1-\alpha _k) F_{I_k}(x^{(k)})+\Big( \bar {\cal A} \Big (\Big (x^{(k+1)}\Big )^{m-1}-\Big (x^{(k)}\Big )^{m-1}
        \Big )\Big )_{I_k} \\
    &\le &(1-\alpha _k) F_{I_k}(x^{(k)})< 0,
\end{eqnarray*}
which shows (\ref{Ik-c}).
For any $i\not\in I_k$, we have $x_i^{(k+1)}=x_i^{(k)}$ and 
\begin{eqnarray*}
F_i(x^{(k+1)}) &=&  a_{i\ldots i}( x_i^{(k+1)} )^{[m-1]} + \sum_{(i_2,\ldots, i_m)\neq (i,\ldots,i)}a_{ii_2\ldots i_m}
    x_{i_2}^{(k+1)}\cdots x_{i_m}^{(k+1)}-b_i\\
     &=&  a_{i\ldots i}( x_i^{(k)} )^{[m-1]} + \sum_{(i_2,\ldots, i_m)\neq (i,\ldots,i)}a_{ii_2\ldots i_m}
    x_{i_2}^{(k+1)}\cdots x_{i_m}^{(k+1)}-b_i\\
    &\le &  a_{i\ldots i}(x_i^{(k)})^{[m-1]}  + \sum_{(i_2,\ldots, i_m)\neq (i,\ldots,i)}a_{ii_2\ldots i_m}
    x_{i_2}^{(k)}\cdots x_{i_m}^{(k)}-b_i \\
    &=& F_i(x^{(k)}).
\end{eqnarray*}
The last inequality yields (\ref{dec-F}).
The proof is complete.
\qed

We conclude this section by showing that the limit of $\{x^{(k)}\}$ exists and
is a solution to the M-TCP (\ref{M-tcp}).
\begin{theorem}\label{th:conv}
Suppose that $\cal A$ is a strong M-tensor. Then the sequence of iterates $\{x^{(k)}\}$ generated by
Algorithm 1 converges to a solution to the M-TCP  (\ref{M-tcp}).
\end{theorem}
{\bf{Proof}}
It follows from the last proposition that there is an index $\bar k$ such that
\[
I_k=I_{\bar k}\stackrel\triangle {=}\bar I,\qquad\forall k\ge \bar k.
\]
By the definition of $I_k$, we always have for $k\ge \bar k$,
\[
x^{(k)}_{\bar I} \in {\cal F}_{\bar I}\stackrel\triangle{=} \{x\ge 0\;| {\cal A}_{\bar I}x_ {\bar I}^{m-1}-b_{\bar I}\le 0\}.
\]
Since $\cal A$ is a strong M-tensor, so is its principal subtensor ${\cal A}_{\bar I}$. It follows from
Lemma \ref{lem-g} that the set ${\cal F}_{\bar I}$ has a largest element $\bar x_{\bar I}\ge x_{I_k}^{(k)}$.
It together with Proposition \ref{prop:monotone-1-c} shows that the
sequence of iterates $\{x_{\bar I}^{(k)}\}$ is monotone and bounded
and hence has a limit $\bar x$.
Taking limits in both sides of (\ref{sub:Meq}), we get
$F_{\bar I}(\bar x_{\bar I})=0$.
On the other hand, by the definition of $x^{(k)}$ and $I_k$,
we always have $x_i^{(k)}=0$ and $F_i( x^{(k)})\ge 0$, $\forall i\not\in \bar I$.
Consequently, we claim that $\bar x=(\bar x_{\bar I}, 0)$ is a solution to the TCP (\ref{M-tcp}).
\qed

\section{An Improvement}\label{sec:3}
In this section, we make some improvements to the method
proposed in the last section. The idea is similar to the approximate Newton method
for solving M-tensor equation by Li, Guan and Wang \cite {Li-Guan-Wang-18}.

Instead of solving the lower dimensional linear system (\ref{sub:Meq}), we
solve the following lower dimensional system of linear equations
\begin{equation}\label{Meq-improved}
  A_{I_k}\Big(x_{I_k}^{(k+1)}\Big)^{[m-1]}-A_{I_k}\Big (x_{I_k}^{(k)}\Big ) ^{[m-1]}+\delta_{I_k} =0,\ \ k = 0,1,2,\ldots ,
\end{equation}
where
\[
\delta_{I_k} = \alpha_kF_{I_k}(x^{(k)})+\epsilon_{I_k}, \quad \alpha_k\in (0,1)
\]
and $\epsilon_{I_k}$ is chosen in the way that $\epsilon_{I_0}=0$ and for $k\geq 1$,
\[
\epsilon_{I_k}=r_{I_k}(x^{(k)})-r_{I_k}(x^{(k-1)}),
\]
where
\[
r_{I_k}(x^{(k)}) = \frac{1}{m-1}(F_{I_k}(x^{(k)})+b_{I_k})-A_{I_k}(x^{(k)})^{[m-1]}
\]
The difference between (\ref{sub:Meq}) and (\ref{Meq-improved})  lies in the term
$\epsilon_{I_k}=r_{I_k}(x^{(k)})-r_{I_k}(x^{(k-1)})$. The role of this term is to
make the solution of (\ref{Meq-improved}) closer to the point generated by Newton's method
for solving $F_{I_k}(x)=0$ in $x_{I_k}$ than the solution of
(\ref{sub:Meq}). We refer to \cite{Li-Guan-Wang-18} for details.

In order for the method to be monotonically convergent,  we need the requirement
$\alpha_kF_{I_k}(x^{(k)})+\epsilon_{I_k}\leq 0$, $\forall k$ to ensure the monotone property of $\{x_{I_k}^{(k)}\}$.
It is satisfied if we let $\epsilon_{I_k}$ satisfy
\begin{equation}\label{eps-improved-add}
\epsilon_{I_k} \leq -\alpha_kF_{I_k}(x^{(k)})\stackrel\triangle {=}\epsilon_{I_k}^+.
\end{equation}
On the other hand, we also need the condition $(1-\alpha _k) F_{I_k}(x^{(k)})-\epsilon_{I_k}\leq 0$.
It is satisfied if we let
 $\epsilon_{I_k}$ satisfy
\begin{equation}\label{eps-improved-minus}
\epsilon_{I_k} \geq (1 -\alpha_k)F_{I_k}(x^{(k)})\stackrel\triangle {=}\epsilon_{I_k}^-.
\end{equation}
Combine (\ref{eps-improved-add}) and (\ref{eps-improved-minus}), we need to choose $\epsilon_{I_k}$
that satisfies both of the inequalities above, i.e.,
\[
\epsilon_{I_k}^- \leq \epsilon_{I_k} \leq \epsilon_{I_k}^+.
\]
Practically, we can first let
$\tilde \epsilon_{I_k} = \min\{\epsilon_{I_k}^+,
r_{I_k}(x_{I_k}^{(k)})-r_{I_k}(x_{I_k}^{(k-1)})\}$ and
let $\epsilon_{I_k} = \max\{\tilde \epsilon_{I_k},\epsilon_{I_k}^-\}$.
And then solve the system of linear equation (\ref{Meq-improved}) to
get a $\bar{x}_{I_k}^{(k+1)}$. Let $\bar{x}_i^{(k+1)}=0$, $\forall i\not\in I_k$.
 If $F_{I_k}(\bar {x}^{(k+1)})<  0$,
then we let $x_{I_k}^{(k+1)}=\bar{x}_{I_k}^{(k+1)}$ and
update the index set ${I_{k}}$ to get ${I_{k+1}}=\{i\;|\;F_i(x^{(k+1)})<0\}$.
Otherwise, we let $\epsilon_{I_k}=0$ and solve (\ref{Meq-improved}) to get a $x_{I_k}^{(k+1)}$.

Based on the above arguments, we propose a lower dimensional approximate Newton approach as follows.\\
\noindent
{\bf{Algorithm 2. (A Sequential Lower Dimensional Approximate Newton Approach)}}
\begin{itemize}
\item [] {\bf Initial.} Given positive sequence $\{\alpha_k\}\subset (\alpha_{\min},1)$
and positive constant $\eta > 0 $.
Find an initial point  $x^{(0)}\ge 0$ such that
the set $I_0=\{i\;|\;F_i(x^{(0)})<0\}$ is not empty.
Let $\epsilon_{I_0}=0$ and $k:=0$.
\item [] {\bf Step 1.} Stop if $\|\min\{F(x^{(k)}),x^{(k)}\}\|\leq \eta$.
\item [] {\bf Step 2.} Solve the lower dimensional system of linear equations (\ref{Meq-improved})
to get $\bar{x}^{(k+1)}_{I_k}$.  Let $\bar{x}_i^{(k+1)}=0$, $\forall i\not\in I_k$.
\item [] {\bf Step 3.} If $F_{I_k}(\bar{x}^{(k+1)})< 0$, then go to Step 4. Otherwise, let
$\epsilon_{I_k}=0$. Go to Step 2.
\item [] {\bf Step 4.} Let $x^{(k+1)}_{I_k}=\bar{x}^{(k+1)}_{I_k}$ and update the index
set $I_{k}$ to get $I_{k+1}=\{i\;|\; F_i(x^{(k+1)})<0\}$. Set $k:=k+1$, and  let
$\tilde\epsilon_{I_k} = \min\{\epsilon_{I_k}^+,
r_{I_k}(x_{I_k}^{(k)})-r_{I_k}(x_{I_k}^{(k-1)})\}$
and $\epsilon_{I_k} = \max\{\tilde \epsilon_{I_k},\epsilon_{I_k}^-\}$.  Go to Step 1.
\end{itemize}

The approximate Newton method possesses
similar properties to the lower dimensional linear equation method proposed
in the last section. In particular, the sequence of iterates generated by the
method converges to a solution to the M-TCP monotonically. We summarize the
results as a theorem below but omit the proof.

\begin{theorem}\label{th:conv-2}
Suppose that ${\cal A}$ is a strong M-tensor. Then Algorithm 2
is well defined. Moreover, the sequence of iterates $\{x^{(k)}\}$ generated
by the algorithm converges to a solution of the M-TCP (\ref{M-tcp}) montonically
in the sense
\begin{equation}\label{monotone-1-improved}
x^{(k+1)}\ge x^{(k)}\ge 0,\quad\forall k\ge 0.
\end{equation}
In addition, the sequence of index sets $\{I_k\}$ is non-decreasing in the sense
\[
I_k\subseteq I_{k+1}.
\]
\end{theorem}

\section{Numerical Results}\label{sec:4}
In this section, we conduct some numerical experiments to test the methods proposed in Sections 2
and 3 on several classes of problems. We implemented our methods in Matlab R2015b and ran the codes on a personal
computer with 3.60 GHz CPU and 20.0 GB RAM.
 We used the tensor toolbox \cite{Bader-Kolda-others-15}
 to proceed the related tensor computation.
The termination criterion for both  methods
  is set to
\[
\|\min\{F(x^{(k)}),x^{(k)}\}\|\leq 10^{-8}
\]
or the number of iteration reaches to 1000. In all cases, we take the parameter $\alpha_k = \alpha$ be constant.

While conducting numerical experiments, we set the initial point to be $x^{(0)} = (0,0,...,0)^T$, and the
 elements of the right hand side $b$ is randomly generated by letting
 $b={\cal A}\tilde x^{m-1}$ with  $\tilde x$ uniformly distributed in $(0,1)$.
The M-tensors of the test problems were  from \cite{Ding-Wei-16,Li-Xie-Xu-17,Xu-Li-Xie-19}.
Details are given below.

\textbf{Problem 1.} The M-tensor ${\cal A}$ takes
 the form ${\cal A}=s{\cal I}-{\cal B}$, where ${\cal B}$ is a  nonnegative $m$-order tensor
  whose elements are random values uniformly  distributed in $(0,1)$. The parameter $s$ is set to
\[
s=(1+\varepsilon)\cdot\max_{i=1,2,...,n}({\cal B}e^{m-1})_i, \ \ \varepsilon>0
\]
where $e=(1,1,...,1)^T$, and we set $\varepsilon=0.01$.
We tested the case $m=3$, $m=4$, and $m=5$ with different dimensions.

\textbf{Problem 2.}  The M-tensor ${\cal A}$ takes
 the form ${\cal A}=s{\cal I}-{\cal B}$, where ${\cal B}$ is a symmetric m-order tensor
whose elements are random values uniformly  distributed in $(0,1)$. The parameter $s$ is set to
\[
s=(1+\varepsilon)\cdot\max_{i=1,2,...,n}({\cal B}e^{m-1})_i, \ \ \varepsilon>0
\]
where $e=(1,1,...,1)^T$, and we set $\varepsilon=0.01$.
We tested the case $m=3$, $m=4$, and $m=5$ with different dimensions.

\textbf{Problem 3.} The M-tensor ${\cal A}$ takes
 the form ${\cal A}=s{\cal I}-{\cal B}$, where ${\cal B}$ is a nonnegative m-order tensor
 whose elements are defined by
\[
b_{i_1i_2...i_m} = |\sin(i_1+i_2+\cdots+i_m)|
\]
and $s=n^{m-1}$. We tested the case $m=3$, $m=4$, and $m=5$ with different dimensions.

We tested both methods with different value $\alpha $
on the above three problems with  $m=3$, $m=4$ and $m=5$ and  different dimensions $n$.
For simplicity, we abbreviate  Algorithms 1 and 2,
as LD-LEQA, and LD-A-Newton respectively.

For each $\alpha$, and each pair $(m,n)$, we tested the methods LD-LEQA and LD-A-Newton 100 times and
recorded the average performance of both methods.
Tables 1 -9 list the performance  of both methods on Problems 1-3 with $m=3$, $m=4$ and $5$
and different sizes $n$. The meaning of each column is given below.

\begin{tabular}{ll}
'$n$': &  the dimension $n$ of the problem;\\
'Iter': &  the total number of iterations; \\
'Time': & the CPU time (in seconds) used for the method;\\
'$I_k$': & the average update times for the index set $I_k$; \\
'Res.': & the average final residual $\|\min\{F(x^{(k)}),x^{(k)}\}\|$;\\
'K': &  the average  cycle times between Steps 2 and 3 in the method LD-A-Newton.
\end{tabular}

First, we point out that the LD-A-Newton method successfully terminated
at a solution of the test problems in all cases, and the LD-LEQA method terminated
at a solution of the test problems in most cases. In the case $m=3, \alpha=0.1$,
the LD-LEQA method failed sometimes. The data in Tables show the
efficiency of the proposed methods.

The convergence  of both  methods  needs the requirement
$\alpha_k<1$.  However, as $\alpha_k$ can be arbitrarily close to 1,
we tested the methods with $\alpha_k=1$. The results show that
$\alpha_k=1$ is practically a good choice.

The results clearly show that both methods with small $\alpha_k$ performed not very well.
In most cases, as the parameter $\alpha_k$ increases, the  performance of both methods becomes better and better.
The best parameter for 'LD-LEQA'method  seems to $\alpha_k=1$ while
the best parameter $\alpha_k$ for the `LD-A-Newton' method seems not to be very clear. In many cases,
$\alpha_k=0.9$ is the best  while for some problems the best $\alpha_k$ seems to $1$. Perhaps the best parameter for that method
lies in the interval $[0.9, 1]$. However, theoretically, it is not easy to
determine the best parameter $\alpha_k$.
The performance of both methods LD-LEQA and LD-A-Newton is very close for the case $\alpha=1$.
However, it does not mean that they coincide in that case.
Compared the best cases for both methods, we can find that in
most cases, the `LD-A-Newton' method is an improvement to the 'LD-LEQA'method though sometimes
the improvement is not too much. At least, the performance of the `LD-A-Newton' method is not at bad as 'LD-LEQA'.

\begin{table}[htbp]
\centering
{\footnotesize
\caption{Results for LD-LEQA and LD-A-Newton on Problem 1 with 3-order tensors}
\begin{tabular}{c|c|cccc|ccccc}
\hline\hline
&  &  & LD-LEQA & & & & & LD-A-Newton & & \\ \hline
n & $\alpha_k$ & Iter. & Time &  ${I_k}$ & Res. & Iter.  & Time & ${I_k}$ & Res. & K\\ \hline
10	&	0.10 	&	411.0 	&	0.00478 	&	1.28 	&	9.83E-09	&	40.7 	&	0.00085 	&	1.27 	&	6.81E-09	&	0	\\
10	&	0.20 	&	207.0 	&	0.00243 	&	1.40 	&	9.49E-09	&	35.1 	&	0.00072 	&	1.35 	&	6.85E-09	&	0	\\
10	&	0.30 	&	142.1 	&	0.00170 	&	1.52 	&	9.27E-09	&	36.4 	&	0.00077 	&	1.46 	&	7.01E-09	&	0	\\
10	&	0.40 	&	102.3 	&	0.00124 	&	1.50 	&	9.04E-09	&	35.2 	&	0.00075 	&	1.40 	&	7.10E-09	&	0	\\
10	&	0.50 	&	76.6 	&	0.00094 	&	1.26 	&	8.34E-09	&	33.0 	&	0.00069 	&	1.23 	&	7.03E-09	&	0	\\
10	&	0.60 	&	67.6 	&	0.00084 	&	1.38 	&	8.49E-09	&	36.4 	&	0.00076 	&	1.31 	&	7.19E-09	&	0	\\
\hline\hline
\end{tabular}}

\end{table}

\begin{table}[htbp]
\centering
{\footnotesize
\caption{Results for LD-LEQA and LD-A-Newton on Problem 1 with 4-order tensors}
\begin{tabular}{c|c|cccc|ccccc}
\hline\hline
&  &  & LD-LEQA & & & & & LD-A-Newton & & \\\hline
n & $\alpha_k$ & Iter. & Time &  ${I_k}$ & Res. & Iter.  & Time &   ${I_k}$ & Res. & K\\\hline
10	&	0.1 	&	307.5 	&	0.00531 	&	0.52 	&	9.65E-09	&	23.6 	&	0.00070 	&	0.49 	&	5.62E-09	&	0	\\
10	&	0.2 	&	150.3 	&	0.00294 	&	0.56 	&	9.30E-09	&	21.2 	&	0.00065 	&	0.56 	&	5.52E-09	&	0	\\
10	&	0.3 	&	91.5 	&	0.00193 	&	0.38 	&	8.87E-09	&	18.6 	&	0.00064 	&	0.38 	&	5.11E-09	&	0	\\
10	&	0.4 	&	70.1 	&	0.00149 	&	0.60 	&	8.38E-09	&	20.3 	&	0.00070 	&	0.56 	&	5.44E-09	&	0	\\
10	&	0.5 	&	51.5 	&	0.00106 	&	0.54 	&	7.97E-09	&	18.7 	&	0.00063 	&	0.54 	&	4.82E-09	&	0	\\
10	&	0.6 	&	41.5 	&	0.00094 	&	0.48 	&	7.60E-09	&	19.2 	&	0.00068 	&	0.48 	&	5.62E-09	&	0	\\
\hline\hline
\end{tabular}}

\end{table}

\begin{table}[htbp]
\centering{\footnotesize
\caption{Results for LD-LEQA and LD-A-Newton on Problem 1 with 5-order tensors}
\begin{tabular}{c|c|cccc|ccccc}\hline\hline
&  &  & LD-LEQA & & & & & LD-A-Newton & & \\ \hline
n & $\alpha_k$ & Iter. & Time &  ${I_k}$ & Res. & Iter.  & Time &  ${I_k}$ & Res. & K\\ \hline
10	&	0.10 	&	289.9 	&	0.00874 	&	0.28 	&	9.61E-09	&	18.9 	&	0.00086 	&	0.27 	&	4.47E-09	&	0	\\
10	&	0.20 	&	136.4 	&	0.00432 	&	0.28 	&	9.15E-09	&	16.0 	&	0.00073 	&	0.27 	&	4.89E-09	&	0	\\
10	&	0.30 	&	86.8 	&	0.00285 	&	0.26 	&	8.64E-09	&	15.3 	&	0.00072 	&	0.26 	&	4.11E-09	&	0	\\
10	&	0.40 	&	63.2 	&	0.00206 	&	0.32 	&	8.11E-09	&	15.3 	&	0.00074 	&	0.31 	&	4.22E-09	&	0	\\
10	&	0.50 	&	47.0 	&	0.00155 	&	0.23 	&	7.53E-09	&	14.8 	&	0.00073 	&	0.23 	&	3.93E-09	&	0	\\
10	&	0.60 	&	37.2 	&	0.00126 	&	0.32 	&	7.05E-09	&	15.1 	&	0.00074 	&	0.32 	&	4.62E-09	&	0	\\
\hline\hline
\end{tabular}}

\end{table}

\begin{table}[htbp]
\centering{\footnotesize
\caption{Results for LD-LEQA and LD-A-Newton on Problem 2 with 3-order tensors}
\begin{tabular}{c|c|cccc|ccccc}\hline\hline
&  &  & LD-LEQA & & & & & LD-A-Newton & & \\ \hline
n & $\alpha_k$ & Iter. & Time &  ${I_k}$ & Res. & Iter.  & Time &  ${I_k}$ & Res. & K\\\hline
10	&	0.10 	&	383.2 	&	0.00437 	&	1.06 	&	8.73E-09	&	37.8 	&	0.00073 	&	1.08 	&	6.70E-09	&	0	\\
10	&	0.20 	&	206.7 	&	0.00242 	&	1.47 	&	9.50E-09	&	35.2 	&	0.00068 	&	1.41 	&	6.96E-09	&	0	\\
10	&	0.30 	&	134.9 	&	0.00161 	&	1.34 	&	9.18E-09	&	34.0 	&	0.00067 	&	1.32 	&	6.92E-09	&	0	\\
10	&	0.40 	&	97.6 	&	0.00118 	&	1.38 	&	8.93E-09	&	33.3 	&	0.00065 	&	1.36 	&	7.02E-09	&	0	\\
10	&	0.50 	&	74.8 	&	0.00093 	&	1.28 	&	8.55E-09	&	32.2 	&	0.00063 	&	1.25 	&	6.87E-09	&	0	\\
10	&	0.60 	&	60.7 	&	0.00076 	&	1.33 	&	8.34E-09	&	32.1 	&	0.00063 	&	1.28 	&	7.05E-09	&	0	\\
\hline\hline
\end{tabular}}

\end{table}

\begin{table}[htbp]
\centering {\footnotesize
\caption{Results for LD-LEQA and LD-A-Newton on Problem 2 with 4-order tensors}
\begin{tabular}{c|c|cccc|ccccc}\hline\hline
&  &  & LD-LEQA & & & & & LD-A-Newton & & \\\hline
n & $\alpha_k$ & Iter. & Time &  ${I_k}$ & Res. & Iter.  & Time &  ${I_k}$ & Res. & K\\ \hline
10	&	0.1 	&	305.9 	&	0.00466 	&	0.43 	&	9.59E-09	&	23.3 	&	0.00061 	&	0.43 	&	4.99E-09	&	0	\\
10	&	0.2 	&	147.3 	&	0.00274 	&	0.52 	&	9.31E-09	&	20.5 	&	0.00059 	&	0.50 	&	5.35E-09	&	0	\\
10	&	0.3 	&	95.0 	&	0.00184 	&	0.52 	&	8.80E-09	&	20.0 	&	0.00061 	&	0.52 	&	5.23E-09	&	0	\\
10	&	0.4 	&	67.3 	&	0.00140 	&	0.51 	&	8.55E-09	&	19.2 	&	0.00062 	&	0.51 	&	5.14E-09	&	0	\\
10	&	0.5 	&	54.3 	&	0.00113 	&	0.65 	&	7.94E-09	&	20.5 	&	0.00067 	&	0.62 	&	5.80E-09	&	0	\\
10	&	0.6 	&	41.1 	&	0.00086 	&	0.50 	&	7.48E-09	&	18.9 	&	0.00065 	&	0.48 	&	5.31E-09	&	0	\\
\hline\hline
\end{tabular}}

\end{table}

\begin{table}[htbp]
\centering {\footnotesize
\caption{Results for LD-LEQA and LD-A-Newton on Problem 2 with 5-order tensors}
\begin{tabular}{c|c|cccc|ccccc} \hline\hline
&  &  & LD-LEQA & & & & & LD-A-Newton & & \\\hline
n & $\alpha_k$ & Iter. & Time &  ${I_k}$ & Res. & Iter.  & Time &  ${I_k}$ & Res. & K\\\hline
10	&	0.10 	&	285.4 	&	0.00622 	&	0.27 	&	9.57E-09	&	18.3 	&	0.00071 	&	0.27 	&	4.16E-09	&	0 	\\
10	&	0.20 	&	140.8 	&	0.00414 	&	0.37 	&	9.06E-09	&	17.3 	&	0.00077 	&	0.35 	&	4.48E-09	&	0 	\\
10	&	0.30 	&	88.5 	&	0.00265 	&	0.33 	&	8.60E-09	&	16.1 	&	0.00074 	&	0.32 	&	4.58E-09	&	0 	\\
10	&	0.40 	&	63.6 	&	0.00193 	&	0.29 	&	8.12E-09	&	15.8 	&	0.00072 	&	0.28 	&	4.72E-09	&	0 	\\
10	&	0.50 	&	47.4 	&	0.00148 	&	0.32 	&	7.74E-09	&	15.1 	&	0.00069 	&	0.31 	&	4.13E-09	&	0 	\\
10	&	0.60 	&	36.9 	&	0.00113 	&	0.27 	&	7.07E-09	&	15.1 	&	0.00068 	&	0.27 	&	4.43E-09	&	0 	\\
\hline\hline
\end{tabular}}

\end{table}

\begin{table}[htbp]
\centering {\footnotesize
\caption{Results for LD-LEQA and LD-A-Newton on Problem 3 with 3-order tensors}
\begin{tabular}{c|c|cccc|ccccc}
\hline\hline
&  &  & LD-LEQA & & & & & LD-A-Newton & & \\\hline
n & $\alpha_k$ & Iter. & Time &  ${I_k}$ & Res. & Iter.  & Time &  ${I_k}$ & Res. & K\\\hline
10	&	0.10 	&	480.3 	&	0.00554 	&	1.94 	&	9.77E-09	&	46.0 	&	0.00097 	&	1.86 	&	7.15E-09	&	0   \\	
10	&	0.20 	&	232.4 	&	0.00271 	&	2.02 	&	9.50E-09	&	40.4 	&	0.00079 	&	1.92 	&	7.09E-09	&	0	\\
10	&	0.30 	&	150.3 	&	0.00179 	&	1.99 	&	9.25E-09	&	38.7 	&	0.00079 	&	1.87 	&	7.56E-09	&	0	\\
10	&	0.40 	&	106.9 	&	0.00130 	&	1.85 	&	8.87E-09	&	37.0 	&	0.00074 	&	1.73 	&	7.04E-09	&	0	\\
10	&	0.50 	&	88.0 	&	0.00108 	&	1.85 	&	8.75E-09	&	38.9 	&	0.00077 	&	1.72 	&	7.38E-09	&	0	\\
10	&	0.60 	&	71.2 	&	0.00088 	&	1.82 	&	8.51E-09	&	38.6 	&	0.00076 	&	1.72 	&	7.32E-09	&	0	\\
\hline\hline
\end{tabular}}

\end{table}

\begin{table}[htbp]
\centering{\footnotesize
\caption{Results for LD-LEQA and LD-A-Newton on Problem 3 with 4-order tensors}
\begin{tabular}{c|c|cccc|ccccc} \hline\hline
&  &  & LD-LEQA & & & & & LD-A-Newton & & \\\hline
n & $\alpha_k$ & Iter. & Time &  ${I_k}$ & Res. & Iter.  & Time &  ${I_k}$ & Res. & K\\\hline
10	&	0.10 	&	336.1 	&	0.00620 	&	0.81 	&	9.68E-09	&	26.4 	&	0.00088 	&	0.80 	&	5.82E-09	&	0	\\										
10	&	0.20 	&	168.4 	&	0.00317 	&	1.02 	&	9.32E-09	&	25.3 	&	0.00078 	&	1.01 	&	5.95E-09	&	0	\\										
10	&	0.30 	&	103.2 	&	0.00197 	&	0.77 	&	8.83E-09	&	22.3 	&	0.00071 	&	0.75 	&	5.62E-09	&	0	\\										
10	&	0.40 	&	76.2 	&	0.00152 	&	0.89 	&	8.57E-09	&	22.9 	&	0.00074 	&	0.86 	&	5.65E-09	&	0	\\										
10	&	0.50 	&	58.2 	&	0.00120 	&	0.89 	&	8.11E-09	&	22.3 	&	0.00075 	&	0.86 	&	5.80E-09	&	0	\\										
10	&	0.60 	&	44.2 	&	0.00090 	&	0.73 	&	7.79E-09	&	20.8 	&	0.00070 	&	0.70 	&	5.40E-09	&	0	\\
\hline\hline
\end{tabular}}

\end{table}

\begin{table}[htbp]
\centering {\footnotesize
\caption{Results for LD-LEQA and LD-A-Newton on Problem 3 with 5-order tensors}
\begin{tabular}{c|c|cccc|ccccc} \hline\hline
&  &  & LD-LEQA & & & & & LD-A-Newton & & \\ \hline
n & $\alpha_k$ & Iter. & Time &   ${I_k}$ & Res. & Iter.  & Time &  ${I_k}$ & Res. & K\\\hline
10	&	0.1 	&	292.8 	&	0.0088 	&	0.33 	&	9.57E-09	&	18.6 	&	0.00084 	&	0.33 	&	4.26E-09	&	0	\\
10	&	0.2 	&	142.1 	&	0.0045 	&	0.37 	&	9.10E-09	&	17.2 	&	0.00080 	&	0.37 	&	4.18E-09	&	0	\\
10	&	0.3 	&	92.2 	&	0.0030 	&	0.46 	&	8.73E-09	&	17.2 	&	0.00082 	&	0.46 	&	4.64E-09	&	0	\\
10	&	0.4 	&	66.2 	&	0.0022 	&	0.39 	&	8.25E-09	&	16.8 	&	0.00082 	&	0.39 	&	4.69E-09	&	0	\\
10	&	0.5 	&	49.3 	&	0.0016 	&	0.40 	&	7.55E-09	&	16.0 	&	0.00078 	&	0.38 	&	4.14E-09	&	0	\\
10	&	0.6 	&	38.6 	&	0.0013 	&	0.38 	&	7.00E-09	&	15.9 	&	0.00077 	&	0.37 	&	4.32E-09	&	0	\\
\hline\hline

\end{tabular}}

\end{table}


\begin{thebibliography}{}

\bibitem {Abraham-94}
Abraham Berman., Robert J. Plemmons.: Nonnegative Matrices in the Mathematical Sciences, classics edition, SIAM, Philadelphia (1994)

\bibitem{Bader-Kolda-others-15}
Bader, B.W., Kolda, T.G., et al.: MATLAB Tensor Toolbox Version 2.6. http://www.sandia. gov/~tgkolda/TensorToolbox/index-2.6. html,(2015)

\bibitem {Bulo-Pelillo-09}
Bul$\grave{o}$, S.R., Pelillo, M.: A game-theoretic approach to hypergraph clustering. In Advances in Neural Information Processing Systems 22. Vancouver, Canada (2009)

\bibitem {Ding-Qi-Wei-13}
Ding, W., Qi, L., Wei, Y.: M-Tensors and Nonsingular M-Tensors. Linear Algebra Appl. 439(10), 3264-3278 (2013)

\bibitem{Ding-Wei-16}
Ding, W., Wei, Y.: Solving multi-linear systems with M-Tensors. J. Sci. Comput. 68(2), 683-715 (2016)

\bibitem {Du-Zhang-19}
Du, S., Zhang, L.: A mixed integer programming approach to the tensor complementarity problem. J. Global Optim. 73(4), 789-800 (2019)

\bibitem{Du-Zhang-Chen-Qi-18}
Du, S., Zhang, L., Chen, C., Qi, L.: Tensor absolute value equations. Sci. China Math. 61(9), 1695-1710 (2018)

\bibitem {Guan-Li-19}
Guan, H.B., Li, D.H.: Linearized Methods for M-Tensor Complementarity Problems, J. Optim. Theory Appl. 184(3), 972-987 (2020)

\bibitem {Han-19}
Han, L.: A continuation method for tensor complementarity problems. J. Optim. Theory Appl. 180(3), 949-963 (2019)

\bibitem {Huang-Qi-17}
Huang, Z.H., Qi, L.: Formulating an n-person noncooperative game as a tensor complementarity problem. Comput. Optim. Appl. 66(3), 557-576 (2017)

\bibitem {Huang-Qi-JOTA-1}
Huang, Z.H., Qi, L.: Tensor complementarity problems-Part I: basic theory. J. Optim. Theory Appl. 183(1), 1-23 (2019)

\bibitem {Huang-Qi-JOTA-2}
Huang, Z.H., Qi, L.: Tensor Complementarity Problems-Part II: Solution Methods. J. Optim. Theory Appl. 183(2), 365-385 (2019)

\bibitem {Huang-Qi-JOTA-3}
Huang, Z.H., Qi, L.: Tensor Complementarity Problems-Part III: Applications. J. Optim. Theory Appl. 183(3), 771-791 (2019)

\bibitem {Li-Guan-Wang-18}
Li, D.H., Guan, H.B., Wang, X.Z.: Finding a Nonnegative Solution to an M-Tensor Equation. (2018) arXiv preprint arXiv:1811.11343

\bibitem {Li-Nie-Zeng-08}
Li, D.H., Nie, Y.Y., Zeng, J.P.: Conjugate Gradient Method for the Linear Complementarity Problem with S-Matrix. Mathematical and Computer Modelling.
48, 918-928 (2008)

\bibitem {Li-Xie-Xu-17}
Li, D.H., Xie, S., Xu, H.R.: Splitting methods for tensor equations. Numer. Linear Algebra Appl. (2017). https://doi.org/10.1002/nla.2102

\bibitem {Liu-Li-Vong-18}
Liu, D., Li, W., Vong, S.W.: Tensor complementarity problems: the GUS-property and an algorithm. Linear Multilinear Algebra. 66(9), 1726-1749 (2018)

\bibitem {Luo-Qi-Xiu-17}
Luo, Z., Qi, L., Xiu, N.: The sparsest solutions to Z-tensor complementarity problems. Optim. Lett. 11(3), 471-482 (2017)

\bibitem {Qi-Chen-Chen-18}
Qi L., Chen H., Chen Y.: Tensor Eigenvalues and their Applications.
vol. 39 of Advances in Mechanics and Mathematics, Springer. Singapore (2018)

\bibitem {Qi-Luo-17}
Qi, L., Luo, Z.: Tensor Analysys, Spectral Theory and Special Tensors, SIAM Philadelphia, 2017.

\bibitem{Song-Qi-16}
Song, Y., Qi, L.: Tensor complementarity problem and semi-positive tensors. J. Optim. Theory Appl. 169(3), 1069-1078 (2016)

\bibitem {Xie-Li-Xu-17}
Xie, S.L., Li, D.H., Xu, H.R.: An iterative method for finding the least solution to the tensor complementarity problem. J. Optim. Theory Appl. 175(1), 119-136 (2017)

\bibitem {Xie-Jin-Wei-2018a}
Xie, Z.J., Jin, X.Q., Wei, Y.M.: Tensor methods for solving symmetric M-tensor systems. J. Sci. Comput. 74(1), 412-425 (2017)

\bibitem {Xu-Li-Xie-19}
Xu, H.R., Li, D.H., Xie, S.L.: An equivalent tensor equation to the tensor complementarity problem with positive semi-definite Z-tensor. Optim. Lett. 13(4), 685-694 (2019)

\bibitem {Zhang-Chen-Zhao-19}
Zhang, K.L., Chen, H.B., Zhao, P.F.: A potential reduction method for tensor complementarity problems. J. Ind. Manag. Optim. 15(2), 429-443 (2019)

\bibitem {Zhao-Fan-18}
Zhao, X., Fan, J.: A semidefinite method for tensor complementarity problems. Optim. Methods Softw. 34(4), 758-769 (2019)

\end{thebibliography}
\end{document}